\documentclass{elsarticle}
\usepackage{amsfonts}
\usepackage{amssymb,amsmath,latexsym}

\def\Box{\vcenter{\vbox{\hrule\hbox{\vrule
     \vbox to 8.8pt{\hbox to 10pt{}\vfill}\vrule}\hrule}}}

\newtheorem{thm}{Theorem}[section]
\newtheorem{lem}[thm]{Lemma}
\newtheorem{cor}[thm]{Corollary}

\newtheorem{prop}[thm]{Proposition}

\newtheorem{remark}{Remark}

\def\endproofbox{\hskip 1.1em\hfill\rule{3pt}{9pt}}
 \newenvironment{proof}%
 {%
 \noindent{\bf Proof:}
 }%
 {%
 \quad\hfill\endproofbox\vspace*{2ex}
 }

\begin{document}

\begin{frontmatter}

\title{Complete permutation polynomials over finite fields of odd characteristic}

\author{Guangkui Xu$^{a,b}$, Xiwang Cao$^a$, Ziran Tu$^c$, Xiangyong Zeng$^d$, Lei Hu$^e$  }

\address{
\noindent$^a$School of Mathematical Sciences, Nanjing University of Aeronautics and Astronautics, Nanjing 210016, China\\
\noindent$^b$Department of Mathematics and Computational Science, Huainan Normal University, Huainan 232038, China\\
\noindent$^c$School of Mathematics and Statistics, Henan University of Science and Technology, Luoyang 471003, China\\
\noindent$^d$Faculty of Mathematics and Statistics, Hubei University, Wuhan 430062, China\\
\noindent$^e$State Key Laboratory of Information Security, Institute of Information Engineering, Chinese Academy of Sciences, Beijing 100093, China
}
\ead{xuguangkuis@gmail.com, xwcao@nuaa.edu.cn, naturetu@gmail.com, xyzeng@hubu.edu.cn, hu@is.ac.cn}
%\thanks{}
%\thanks{The work of this paper was supported by National Natural Science Foundation of China (Grant Number: 10971250)}
\date{}
%\maketitle
\begin{abstract}
In this paper, we present three classes of complete permutation monomials over finite fields
  of odd characteristic. Meanwhile, the compositional inverses of these complete permutation polynomials are also proposed.
\end{abstract}


\end{frontmatter}

%\maketitle

\section{Introduction}
\label{intro}

 Let $p$ be a prime number and $q=p^n$. Let $\mathbb{F}_{q}$ denote the finite field of order $q$ and $\mathbb{F}_{q}^{\ast}$ the set of  all non-zero elements of $\mathbb{F}_{q}$.  A
 polynomial $f(x)\in \mathbb{F}_{q} [x]$ is called a {\it permutation polynomial} (PP) of $\mathbb{F}_{q}$ if the associated polynomial function $f: c \mapsto f(c)$ from $\mathbb{F}_{q}$ to $\mathbb{F}_{q}$ is a permutation of $\mathbb{F}_{q}$. For a permutation  polynomial $f(x)\in \mathbb{F}_{q} [x]$ there exists  (a unique) $ f^{-1}(x)\in \mathbb{F}_{q} [x]$ such that $f(f^{-1}(x))\equiv f^{-1}(f(x))\equiv x ({\rm mod}\  x^q-x)$. We call $f^{-1}(x)$ the {\it compositional inverse} of $f(x)$. Permutation polynomials were studied first by Hermite \cite{Hermite} and later by Dickson \cite{Dickson}. Permutation polynomials have been an active topic of study in recent years due to  their important applications in cryptography, coding theory, combinatorial designs  theory. A permutation polynomial $f(x) \in \mathbb{F}_{q}[x]$  is  a {\it complete permutation polynomial} (CPP) over $\mathbb{F}_{q}$ if $f(x)+x$ permutes $\mathbb{F}_{q}$ as well. The study of complete permutation
polynomials started with the work of Niederreiter and Robinson \cite{Niederreite}. Finding new PPs and CPPs of finite fields is a difficult problem and there are rare classes of CPPs known. More investigations on PPs and CPPs  can be found in \cite{Akbary1,Akbary2, Cao,Charpin,Ding, Hou,Laigle,Payne,Zeng, Yuan,Zieve}.

Our interest in complete permutation polynomials arises from a recent paper
 by Tu et al.  \cite{Tu} in which several classes of complete permutation polynomials over finite fields
of even characteristic were constructed. More precisely,
they considered  three classes of complete permutation monomials and a class of
trinomial complete permutation polynomials. In \cite{Niederreite}, Niederreiter and Robinson pointed out that the compositional inverse of a complete permutation polynomial is also a complete permutation polynomial.  As  one of our main results in this paper we present three
new classes of monomial complete permutations over finite fields of odd characteristic, not corresponding to any known monomial complete permutation. In order to prove the complete permutation behavior of the second class of monomials, some properties of Dickson polynomials will be employed.

The rest of this paper is organized as follows.  Some preliminaries
and notations are given in Section 2. In Section 3,  we propose three classes of monomial complete permutations over finite fields of odd characteristic and present their compositional inverses. %Conclusions are given in Section 4.

\section{Notations and preliminaries }
\label{sec:1}
 Let $p$ be a prime number and $q=p^n$. For any positive  integer $n$ with a divisor $m\geq 1$, the
trace function, denoted by ${\rm Tr}_m^n(x)$, from $\mathbb{F}_{p^{n}}$ to
$\mathbb{F}_{p^{m}}$ is defined as
$${\rm Tr}_{m}^{n}(x)=x+x^{p^{m}}+x^{p^{2m}}+\cdots+x^{p^{(n/m-1)m}}.$$

The determination of permutation polynomials is a nontrivial problem, and some simple examples of permutation polynomials can be obtained from the following result.
{\lem \label{lem:2.1} {\rm\cite{Lidl}} The monomial $x^n$ is a permutation polynomial of $\mathbb{F}_{q}$ if and only if $\gcd(n,q-1)=1$.}

A well-known criterion for  permutation polynomial which will be frequently  used in this paper is the
following lemma:
{\lem \label{lem:2.2}{\rm\cite{Lidl}} The polynomial $f\in\mathbb{F}_{q}[x]$ is a  permutation polynomial of $\mathbb{F}_{q}$ if and only if for every nonzero $\gamma \in\mathbb{F}_{q}$,
\begin{equation}
\sum_{x\in \mathbb{F}_{q}}\omega^{{\rm Tr}^{n}_{1}(\gamma f(x))}=0,
\end{equation}
where $\omega$ is a primitive $p$-th root of unity.}

 Now we recall the knowledge of Dickson polynomials over $\mathbb{F}_{q}$. Dickson polynomials
are a special source of permutation polynomials over finite fields. The reader can refer to the monograph of Lidl, Mullen and Turnwald \cite{Mullen}  for many useful properties and applications of Dickson polynomials. A Dickson polynomial is defined by
 \begin{equation*}
D_{n}(x,a)=\sum_{i=0}^{\lfloor\frac{n}{2}\rfloor}\frac{n-i}{n}{\ n-i \choose i} (-a)^{i}x^{n-2i},
\end{equation*}
where $a\in \mathbb{F}_{q}$, $\lfloor n/2\rfloor$ is the floor function, i.e., the biggest integer less than or equal to $n/2$, and ${\ n-i \choose i}$ is the combinatorial number of $n-i$ chooses $i$.

Further, the family of Dickson polynomials $D_{n}(x,a)\in \mathbb{F}_{q}[x]$  can also be defined by the recurrence
relation
$$D_{i+2}(x,a)=xD_{i+1}(x,a)-aD_{i}(x,a),i=0,1,\cdots$$
with initial values
$$D_{0}(x,a) = 2, D_{1}(x,a)=x.$$
For example, the first few Dickson polynomials over $\mathbb{F}_{3^{m}}$  are given below.
 \begin{eqnarray*}
&&D_{2}(x,a)=x^{2}-2a,
\\&&D_{3}(x,a)=x^{3},
\\&&D_{4}(x,a)=x^{4}-ax^{2}+2a^{2}x,
\\&&D_{5}(x,a)=x^{5}+ax^{3}-a^{2}x.
\end{eqnarray*}
We also have the following fundamental result.
{\lem \label{lem:2.3}{\rm\cite{Lidl}} For a nonzero element $a\in \mathbb{F}_{q}$, Dickson polynomial $D_{n}(x,a)$ over $\mathbb{F}_{q}$ is a permutation polynomial if and only if $\gcd(n,q^{2}-1)=1$.}

In Section 3, we will show that the complete permutation polynomials in the second class
 are related to some properties of
Dickson polynomials.

 The following two lemmas will be used in Section 3.
{\lem \label{lem:2.5}{\rm\cite{Niederreite}} Let $f(x)$ be a complete permutation polynomial over $\mathbb{F}_{q}$. Then $f^{-1}(x)$ is also a complete permutation polynomial over $\mathbb{F}_{q}$.}

{\lem \label{lem:2.4}{\rm\cite{Wan}} Pick $d>0$ with $d|q-1$, and let $\zeta$ be a primitive $d$-th root of unity in $\mathbb{F}_{q}$. Then the polynomial $x^{\frac{q-1}{d}+1}+ax(a\neq0)$ is a permutation polynomial of $\mathbb{F}_{q}$ if and only if the following conditions are satisfied:

\noindent(i) $(-a)^{d}\neq1$;

\noindent(ii) For all $0 \leq i < j \leq d-1$,
\begin{eqnarray*}
\left(  \frac{a+\zeta^{i}}{a+\zeta^{j}}\right)^{\frac{q-1}{d}}\neq \zeta^{j-i}.
\end{eqnarray*}}

\section{Three classes of monomial CPPs   over finite fields of odd characteristic  }
In this section, three classes of monomial polynomials over finite fields of odd characteristic  are explored. The study of these monomials will start with a technique used
by Dobbertin \cite{Dobbertin}, Leander \cite{Leander} and Tu et al.\cite{Tu}.

 We fix $p$ as an odd prime number in this section. Let the integer $n=2m$ for an odd integer $m$. Since $m$ is odd, the polynomial $x^2+1$ is irreducible over $\mathbb{F}_{3^m}$ as it is irreducible over $\mathbb{F}_{3}$. Let $\alpha$ be a root of  $x^2+1$. Then the order of $\alpha$ is $4$ in the multiplicative group of  $\mathbb{F}_{3^{2m}}=\mathbb{F}_{3^m}(\alpha)$. In the sequel let
$$x=x_{0}+x_{1}\alpha, \ \ x_{0},x_{1}\in \mathbb{F}_{3^m}$$
be an arbitrary element of $\mathbb{F}_{3^{2m}}$. Since $m$ is odd, we have $3^m\equiv 3 ({\rm mod}\  4)$ and $\alpha^{3^{m}}=\alpha^{3}$. We conclude that
\begin{eqnarray}
{\rm Tr}_{m}^{2m}(\alpha)={\rm Tr}_{m}^{2m}(\alpha^{3})=0, \ \ {\rm Tr}_{m}^{2m}(\alpha^{2})=1,
\end{eqnarray}
and therefore
\begin{eqnarray}
{\rm Tr}_{m}^{2m}(x)={\rm Tr}_{m}^{2m}(x_{0}+x_{1}\alpha)=2x_{0}.
\end{eqnarray}
{\thm\label{thm:3.1}For any positive odd  integer $m$ and a nonzero element $v$ in $\mathbb{F}_{3^{2m}}$ with ${\rm Tr}_{m}^{2m}(v)=0$, the monomial $v^{-1}x^{3^{m}+2}$ is a complete permutation polynomial over $\mathbb{F}_{3^{2m}}$.}

\begin{proof}
 Denote
\begin{equation}\label{f-11121}
   S=\{v_{0}+v_{1 }\alpha: v_{0},v_{1 }\in \mathbb{F}_{3^{m}},v_{0}=0 \}\backslash\{0\},
\end{equation}
where $\alpha$ is defined as above. By Eqs.(2) and (3), we know that $S$ is the set of all nonzero elements $v$ in $\mathbb{F}_{3^{2m}}$ with ${\rm Tr}_{m}^{2m}(v)=0$. For each $v\in S $, from Lemma \ref{lem:2.1}, the monomial $v^{-1}x^{3^{m}+2}$ is a permutation polynomial over $\mathbb{F}_{3^{2m}}$, since $\gcd(3^{2m}-1,3^{m}+2)=\gcd(3^m-1,3^{m}+2)=\gcd(3^m-1,3)=1$. To prove $v^{-1}x^{3^{m}+2}$ is a CPP over $\mathbb{F}_{3^{2m}}$, it is sufficient to show that $x^{3^{m}+2}+vx$ is a PP over $\mathbb{F}_{3^{2m}}$ for each $v\in S $.

 Note that $\gcd(3^{2m}-1,3^{m}+2)=1$, hereafter the nonzero $\gamma \in \mathbb{F}_{3^{2m}}$ will be represented as $\gamma=\beta^{3^{m}+2}$ for a unique nonzero $\beta \in \mathbb{F}_{3^{2m}}$. Then we have
\begin{eqnarray*}
&&\sum_{x\in \mathbb{F}_{3^{2m}}}\omega^{{\rm Tr}_{1}^{2m}\left(\gamma(x^{3^{m}+2}+vx)\right)}\\
&=&\sum_{x\in \mathbb{F}_{3^{2m}}}\omega^{{\rm Tr}_{1}^{2m}\left((\beta x)^{3^{m}+2}+\beta^{3^{m}+1}v(\beta x)\right)}
\\&=&\sum_{x\in \mathbb{F}_{3^{2m}}}\omega^{{\rm Tr}_{1}^{2m}( x^{3^{m}+2}+\beta^{3^{m}+1}vx)}
\\&=&\sum_{x\in \mathbb{F}_{3^{2m}}}\omega^{{\rm Tr}_{1}^{m}\left({\rm Tr}_{m}^{2m}( x^{3^{m}+2}+\beta^{3^{m}+1}vx)\right)}.\end{eqnarray*}
By expressing $x\in \mathbb{F}_{3^{2m}}$ as $x_{0}+x_{1}\alpha$ and Eq.(2),  we compute
\begin{eqnarray*}
&&{\rm Tr}_{m}^{2m}( x^{3^{m}+2})\\&=&{\rm Tr}_{m}^{2m}\left( (x_{0}+x_{1}\alpha)^{3^{m}+2}\right)
\\&=&{\rm Tr}_{m}^{2m}\left( (x_{0}+x_{1}\alpha^{3^{m}})(x_{0}+x_{1}\alpha)^{2}\right)
\\&=&{\rm Tr}_{m}^{2m}\left( x_{0}^{3}+2x_{0}x_{1}^{2}+(2x_{0}^{2}x_{1}+x_{1}^{3})\alpha+ x_{0}x_{1}^{2}\alpha^{2}+ x_{0}^{2}x_{1}\alpha^{3}\right)
\\&=&2(x_{0}^{3}+x_{0}x_{1}^{2})
\end{eqnarray*}
since $\alpha^{3^{m}}=\alpha^{3}$ for odd $m$. %Applying the property ${\rm Tr}_{1}^{m}(x)={\rm Tr}_{1}^{m}(x^3)$, one has that
%\begin{eqnarray*}
%\sum_{x\in \mathbb{F}_{3^{2m}}}\omega^{{\rm Tr}_{1}^{2m}(x^{3^{m}+2})}&=&\sum_{x\in \mathbb{F}_{3^{2m}}}\omega^{{\rm Tr}_{1}^{m}({\rm Tr}_{m}^{2m}( x^{3^{m}+2}))}
%\\&=&\sum_{x_{0},x_{1}\in \mathbb{F}_{3^{m}}}\omega^{{\rm Tr}_{1}^{m}(2(x_{0}^{3}+x_{0}x_{1}^{2}))}
%\\&=&\sum_{x_{0},x_{1}\in \mathbb{F}_{3^{m}}}\omega^{{\rm Tr}_{1}^{m}(2x_{0}^{3}(x_{1}^{6}+1))}.
%\end{eqnarray*}
%Since $x^{3^{m}+2}$ is a PP over $\mathbb{F}_{3^{2m}}$, then $\sum\limits_{x\in \mathbb{F}_{3^{2m}}}\omega^{{\rm Tr}_{1}^{2m}(x^{3^{m}+2})}=0$ by Lemma \ref{lem:2.2}.

Note that $(\beta^{3^{m}+1})^{3^{m}-1}=1$, we have $\beta^{3^{m}+1}\in\mathbb{F}_{3^{m}}$ and $\beta^{3^{m}+1}v \in S$. By Eq. (\ref{f-11121}), we can assume that $\beta^{3^{m}+1}v=u=u_{1}\alpha$ with $u_{1}\in \mathbb{F}_{3^{m}}$, and then
\begin{eqnarray*}
{\rm Tr}_{m}^{2m}( \beta^{3^{m}+1}vx)
={\rm Tr}_{m}^{2m}\left( u_{1}\alpha(x_{0}+x_{1}\alpha)\right)
={\rm Tr}_{m}^{2m}(u_{1}x_{0}\alpha+u_{1}x_{1}\alpha^{2})
=u_{1}x_{1}.
\end{eqnarray*}
Combining with the fact ${\rm Tr}_{1}^{m}(z^{3})={\rm Tr}_{1}^{m}(z)$ for any $z\in \mathbb{F}_{3^{m}}$, we have
\begin{eqnarray*}
&&\sum_{x\in \mathbb{F}_{3^{2m}}}\omega^{{\rm Tr}_{1}^{2m}\left(\gamma (x^{3^{m}+2}+vx)\right)}\\&=&\sum_{x\in \mathbb{F}_{3^{2m}}}\omega^{{\rm Tr}_{1}^{m}\left({\rm Tr}_{m}^{2m}( x^{3^{m}+2}+\beta^{3^{m}+1}vx)\right)}
\\&=&\sum_{x_{0},x_{1}\in \mathbb{F}_{3^{m}}}\omega^{{\rm Tr}_{1}^{m}\left(2x_{0}^{3}(x_{1}^{6}+1)+ u_{1}x_{1} \right)}
\\&=&\sum_{x_{1}\in \mathbb{F}_{3^{m}}}\omega^{{\rm Tr}_{1}^{m}(u_{1}x_{1})}\sum_{x_{0}\in \mathbb{F}_{3^{m}}}\omega^{{\rm Tr}_{1}^{m}\left(2x_{0}^{3}(x_{1}^{6}+1)\right )}
\\&=&0
\end{eqnarray*}
since the equation $ x_1^{6}+1=0$ has no solution in $\mathbb{F}_{3^{m}}$ ($-1$ is a non-square element in $\mathbb{F}_{3^{m}}$ for odd $m$).

Hence, for every nonzero $\gamma \in\mathbb{F}_{3^{2m}}$, we have $$\sum \limits_{x\in \mathbb{F}_{3^{2m}}}\omega^{{\rm Tr}_{1}^{2m}\left(\gamma (x^{3^{m}+2}+vx)\right)}=0.$$
 By Lemma \ref{lem:2.2}, the assertion is proved.
\end{proof}
{\remark \label{remark:1} Besides $x^2+1$, we can also use the other two irreducible polynomials of degree 2 over $\mathbb{F}_{3}$ to prove Theorem \ref{thm:3.1}.
In the case of $x^2+2x+2$, the corresponding set $S$ is  $ S=\{v_{0}+v_{1 }\alpha: v_{0},v_{1 }\in \mathbb{F}_{3^{m}},v_{0}=v_1 \}\backslash\{0\}$. In the case of $x^2+x+2$, the related set $S$ should be $ S=\{v_{0}+v_{1 }\alpha: v_{0},v_{1 }\in \mathbb{F}_{3^{m}},v_1 =2v_{0} \}\backslash\{0\}$.}

{\prop \label{prop:3.2}For any positive odd  integer $m$ and $v$ in $\mathbb{F}_{3^{2m}}^*$ with ${\rm Tr}_{m}^{2m}(v)=0$, the monomial $v^{2\cdot3^{2m-1}-3^{m-1}}x^{2\cdot3^{2m-1}-3^{m-1}}$ is a complete permutation polynomial over $\mathbb{F}_{3^{2m}}$.}

 \begin{proof}
 In Lemma \ref{lem:2.5}, put $f(x)=v^{-1}x^{3^{m}+2}$.
 Observe that
 \begin{eqnarray*}
   && (3^{m}+2)(2\cdot3^{2m-1}-3^{m-1})\\&=&2\cdot3^{3m-1}+3^{2m }-2\cdot3^{m-1}
\\ &\equiv&  1 ({\rm mod}\  3^{2m}-1).
 \end{eqnarray*}
 We get that

 $$f^{-1}(x)=v^{2\cdot3^{2m-1}-3^{m-1}}x^{2\cdot3^{2m-1}-3^{m-1}} .$$
 This leads to the claimed result  from Lemma \ref{lem:2.5}.
 \end{proof}

 In what follows, we will propose  the second class of complete permutation monomials  over finite fields of characteristic 3 based on the properties of Dickson polynomials.

First, we start with a similar analysis as in Theorem \ref{thm:3.1}. Let $n= 2m$, where $m$ is odd. Then $\gcd(3^{2m}-1,2\cdot3^{m}+3)=\gcd(3^m-1,2\cdot3^{m}+3)=\gcd(3^m-1,5)=1 $.  Since $m$ is odd, the polynomial $x^{2}+2x+2$  is irreducible over $\mathbb{F}_{3^{m}}$ as it is irreducible over $\mathbb{F}_{3}$. Let $\alpha$  be a root of $x^{2}+2x+2$. Then  $\alpha$ is a primitive element of $\mathbb{F}_{9}$ and
$\mathbb{F}_{3^{2m}}=\mathbb{F}_{3^{m}}(\alpha)$. Thus, each $x\in \mathbb{F}_{3^{2m}}$ can be written as
\begin{eqnarray*}
x=x_{0}+x_{1}\alpha,\ \ x_{0}, x_{1} \in \mathbb{F}_{3^{m}}.
\end{eqnarray*}
Because $m$ is odd, $3^{m}\equiv 3 ({\rm mod}\ 8)$, and then $\alpha^{3^{m}}=\alpha^{3}$.  We have
\begin{eqnarray}
{\rm Tr}_{m}^{2m}(\alpha)={\rm Tr}_{m}^{2m}(\alpha^{3})=1, {\rm Tr}_{m}^{2m}(\alpha^{2})={\rm Tr}_{m}^{2m}(\alpha^{6})=0,
\end{eqnarray}
 and therefore
\begin{eqnarray}
{\rm Tr}_{m}^{2m}(x)={\rm Tr}_{m}^{2m}(x_{0}+x_{1}\alpha)=2x_{0}+x_{1}.
\end{eqnarray}
{\thm \label{thm:3.2}For any positive odd  integer $m$, the monomial $v^{-1}x^{2\cdot3^{m}+3}$ is a complete permutation polynomial over $\mathbb{F}_{3^{2m}}$ if $v$ is a nonzero element in $\mathbb{F}_{3^{2m}}$ with ${\rm Tr}_{m}^{2m}(\alpha v)=0$ or ${\rm Tr}_{m}^{2m}(\alpha^{3} v)=0$, where $\alpha\in \mathbb{F}_{3^{2m}}$ is a root of the equation $x^2+2x+2=0$.}

\begin{proof}
 Denote
$$ S=\{v_{0}+v_{1 }\alpha: v_{0},v_{1 }\in \mathbb{F}_{3^{m}},v_{0}=0\  {\rm or}  \ v_{1}=2v_{0} \}\backslash\{0\},$$
where $\alpha$ is a root of $x^{2}+2x+2$. By Eqs.(5) and (6), we conclude that $S$ is the set of all nonzero elements $v$ in $\mathbb{F}_{3^{2m}}$ with ${\rm Tr}_{m}^{2m}(\alpha v)=0$ or ${\rm Tr}_{m}^{2m}(\alpha^{3} v)=0$. Note that $\gcd(3^{2m}-1,2\cdot3^{m}+3)=1$, the monomial $v^{-1}x^{3^{m}+2}$ is a permutation polynomial over $\mathbb{F}_{3^{2m}}$ by Lemma \ref{lem:2.1} for every $v\in S $.  The rest of the proof is to show that $x^{2\cdot3^{m}+3}+vx$ permutes $\mathbb{F}_{3^{2m}}$ for every $v\in S $.

  Note that $\gcd(3^{2m}-1,2\cdot3^{m}+3)=1$. Then each $\gamma \in \mathbb{F}_{3^{2m}}^*$ is uniquely written as $\beta^{2\cdot3^{m}+3}$ for a $\beta \in \mathbb{F}_{3^{2m}}^*$. We have
\begin{eqnarray*}
&&\sum_{x\in \mathbb{F}_{3^{2m}}}\omega^{{\rm Tr}_{1}^{2m}\left(\gamma(x^{2\cdot3^{m}+3}+vx)\right)}\\
&=&\sum_{x\in \mathbb{F}_{3^{2m}}}\omega^{{\rm Tr}_{1}^{2m}\left((\beta x)^{2\cdot3^{m}+3}+\beta^{2\cdot3^{m}+2}v(\beta x)\right)}
\\&=&\sum_{x\in \mathbb{F}_{3^{2m}}}\omega^{{\rm Tr}_{1}^{2m}( x^{2\cdot3^{m}+3}+\beta^{2\cdot3^{m}+2}vx)}
\\&=&\sum_{x\in \mathbb{F}_{3^{2m}}}\omega^{{\rm Tr}_{1}^{m}\left({\rm Tr}_{m}^{2m}( x^{2\cdot3^{m}+3}+\beta^{2\cdot3^{m}+2}vx)\right)}.
\end{eqnarray*}
By  Eqs.(5), (6) and $\alpha^{3^{m}} =\alpha^{3}$,  one has that
\begin{eqnarray*}
&&{\rm Tr}_{m}^{2m}( x^{2\cdot3^{m}+3})\\&=&{\rm Tr}_{m}^{2m}\left( (x_{0}+x_{1}\alpha)^{2\cdot3^{m}+3}\right)
\\&=&{\rm Tr}_{m}^{2m}\left( (x_{0}+x_{1}\alpha^{3^{m}})^{2}(x_{0}+x_{1}\alpha)^{3}\right)
\\&=&{\rm Tr}_{m}^{2m}\left( x_{0}^{5}+(2x_{0}^{4}x_{1}+x_{0}^{2}x_{1}^{3})\alpha^{3}+ (2x_{0}x_{1}^{4}+x_{0}^{3}x_{1}^{2})\alpha^{6}+ x_{1}^{5}\alpha\right)
\\&=&2x_{0}^{5}+2x_{0}^{4}x_{1}+x_{0}^{2}x_{1}^{3}+x_{1}^{5}.
\end{eqnarray*}
Because $(\beta^{2\cdot3^{m}+2})^{3^{m}-1}=1$, we have $\beta^{2\cdot3^{m}+2}\in\mathbb{F}_{3^{m}}$ and $\beta^{3^{m}+1}v \in S$. Let $\beta^{2\cdot3^{m}+2}v=u=u_{0}+u_{1}\alpha$ with $u_{0}, u_1\in \mathbb{F}_{3^{m}}$. Then
\begin{eqnarray*}
&&{\rm Tr}_{m}^{2m}( \beta^{3^{m}+1}vx)\\&=&{\rm Tr}_{m}^{2m}( ux)
\\&=&{\rm Tr}_{m}^{2m}\left( (u_{0}+u_{1}\alpha)(x_{0}+x_{1}\alpha)\right)
\\&=&{\rm Tr}_{m}^{2m}\left( u_{0}x_{0}+(u_{0}x_{1}+u_{1}x_{0})\alpha+u_{1}x_{1}\alpha^{2}\right)
\\&=&2u_{0}x_{0}+u_{0}x_{1}+u_{1}x_{0}
\end{eqnarray*}
due to ${\rm Tr}_{m}^{2m}(\alpha)=1$ and ${\rm Tr}_{m}^{2m}(\alpha^{2})=0$.
Thus,
\begin{eqnarray}
&&\sum_{x\in \mathbb{F}_{3^{2m}}}\omega^{{\rm Tr}_{1}^{2m}\left(\gamma(x^{2\cdot3^{m}+3}+vx)\right)}\notag
\\&=&\sum_{x\in \mathbb{F}_{3^{2m}}}\omega^{{\rm Tr}_{1}^{m}\left({\rm Tr}_{m}^{2m}( x^{2\cdot3^{m}+3}+ux)\right)}\notag
\\&=&\sum_{x_{0},x_{1}\in \mathbb{F}_{3^{m}}}\omega^{{\rm Tr}_{1}^{m}(2x_{0}^{5}+2x_{0}^{4}x_{1}+x_{0}^{2}x_{1}^{3}+x_{1}^{5}+2u_{0}x_{0}+u_{0}x_{1}+u_{1}x_{0} )}.
\end{eqnarray}
Case (i). When $u_{0}=0$, Eq.(7) can be rewritten as
\begin{eqnarray}
&&\sum_{x_{0},x_{1}\in \mathbb{F}_{3^{m}}}\omega^{{\rm Tr}_{1}^{m}(2x_{0}^{5}+2x_{0}^{4}x_{1}+x_{0}^{2}x_{1}^{3}+x_{1}^{5}+u_{1}x_{0} )}
\notag\\&=&\sum_{x_{0}\in \mathbb{F}_{3^{m}}}\omega^{{\rm Tr}_{1}^{m}(2x_{0}^{5}+u_{1}x_{0})}\sum_{x_{1}\in \mathbb{F}_{3^{m}}}\omega^{{\rm Tr}_{1}^{m}(x_{1}^{5}+x_{0}^{2}x_{1}^{3}-x_{0}^{4}x_{1} )}
\notag\\&=&\sum_{x_{1}\in \mathbb{F}_{3^{m}}}\omega^{{\rm Tr}_{1}^{m}(x_{1}^{5} )}+\sum_{x_{0}\in \mathbb{F}_{3^{m}}^{*}}\omega^{{\rm Tr}_{1}^{m}(2x_{0}^{5}+u_{1}x_{0})}\sum_{x_{1}\in \mathbb{F}_{3^{m}}}\omega^{{\rm Tr}_{1}^{m}(x_{1}^{5}+x_{0}^{2}x_{1}^{3}-x_{0}^{4}x_{1} )}.
\end{eqnarray}
  Since $\gcd(5,3^{m}-1)=1$,  $x_{1}^{5}$ permutes $\mathbb{F}_{3^{m}}$ by Lemma \ref{lem:2.1} and then
  $$\sum\limits_{x_{1}\in \mathbb{F}_{3^{m}}}\omega^{{\rm Tr}_{1}^{m}(x_{1}^{5} )}=0.$$
Note that polynomial $x_{1}^{5}+x_{0}^{2}x_{1}^{3}-x_{0}^{4}x_{1}$ is a Dickson polynomial of degree 5 over $ \mathbb{F}_{3^{m}}$ for  any nonzero  $x_{0}\in \mathbb{F}_{3^{m}}$. Since $\gcd(5,3^{2m}-1)=1$, by Lemma \ref{lem:2.3},  $x_{1}^{5}+x_{0}^{2}x_{1}^{3}-x_{0}^{4}x_{1}$   permutes $\mathbb{F}_{3^{m}}$ which gives
$$\sum\limits_{x_{1}\in \mathbb{F}_{3^{m}}}\omega^{{\rm Tr}_{1}^{m}(x_{1}^{5}+x_{0}^{2}x_{1}^{3}-x_{0}^{4}x_{1} )}=0$$
for any nonzero  $x_{0}\in \mathbb{F}_{3^{m}}$.
 Consequently, when $u_{0}=0$, the sum in Eq.(8) equals to
   $$\sum_{x_{0},x_{1}\in \mathbb{F}_{3^{m}}}\omega^{{\rm Tr}_{1}^{m}(2x_{0}^{5}+2x_{0}^{4}x_{1}+x_{0}^{2}x_{1}^{3}+x_{1}^{5}+u_{1}x_{0} )}=0. $$
Case (ii). When $u_{1}=2u_{0}$, substituting  $x_{0} $ by $y-x_{1}$ in Eq.(7) yields
\begin{eqnarray*}
&&\sum_{x_{0},x_{1}\in \mathbb{F}_{3^{m}}}\omega^{{\rm Tr}_{1}^{m}\left(2x_0^{5}+2x_0^{4}x_{1}+x_0^{2}x_1^{3}+x_1^{5}+u_{0}(x_0+x_1) \right)}
\notag\\&=&\sum_{y\in \mathbb{F}_{3^{m}}}\omega^{{\rm Tr}_{1}^{m}(2y^{5}+u_0 y)}\sum_{x_{1}\in \mathbb{F}_{3^{m}}}\omega^{{\rm Tr}_{1}^{m}\left(-(x_{1}^{5}+y^{2}x_{1}^{3}-y^{4}x_{1}) \right)}
\notag\\&=&\sum_{x_{1}\in \mathbb{F}_{3^{m}}}\omega^{{\rm Tr}_{1}^{m}(-x_{1}^{5} )}+\sum_{y\in \mathbb{F}_{3^{m}}^{*}}\omega^{{\rm Tr}_{1}^{m}(2y^{5}+u_0 y)}\sum_{x_{1}\in \mathbb{F}_{3^{m}}}\omega^{{\rm Tr}_{1}^{m}\left(-(x_{1}^{5}+y^{2}x_{1}^{3}-y^{4}x_{1})\right )}.
\end{eqnarray*}
Since $x_{1}^{5}+y^{2}x_{1}^{3}-y^{4}x_{1}$ is a Dickosn polynomial of degree $5$ in variable $x_1$ for any fixed $y\in \mathbb{F}_{3^m}$, by a similar analysis as above, we know also that
 $$\sum_{x_{0},x_{1}\in \mathbb{F}_{3^{m}}}\omega^{{\rm Tr}_{1}^{m}\left(2x_0^{5}+2x_0^{4}x_{1}+x_0^{2}x_1^{3}+x_1^{5}+u_{0}(x_0+x_1) \right)}=0 $$
for $u_{1}=2u_{0}$.

Finally, we have $$\sum \limits_{x\in \mathbb{F}_{3^{2m}}}\omega^{{\rm Tr}_{1}^{2m}\left(\gamma (x^{2\cdot3^{m}+3}+vx)\right)}=0$$ for each nonzero $\gamma \in\mathbb{F}_{3^{2m}}$ and so by Lemma \ref{lem:2.2} the desired result is  proved.
\end{proof}
%{\remark \label{remark:3.3} The same result will be obtained if we choose the other primitive polynomial $x^2+x+2$ of degree 2 over $\mathbb{F}_{3}$.}

{\prop \label{prop:3.4}For any positive odd  integer $m$, the monomial $v^{-\frac{2\cdot3^m-3}{5}}x^{-\frac{2\cdot3^m-3}{5}}$ is a complete permutation polynomial over $\mathbb{F}_{3^{2m}}$ for any  nonzero element $v$ in $\mathbb{F}_{3^{2m}}$ with ${\rm Tr}_{m}^{2m}(\alpha v)=0$ or ${\rm Tr}_{m}^{2m}(\alpha^{3} v)=0$, where $\alpha\in \mathbb{F}_{3^{2m}}$ is a root of the equation $x^2+2x+2=0$.}

\begin{proof}
Set $f(x)=v^{-1}x^{2\cdot3^{m}+3}$.  Note that
\begin{eqnarray*}
(2\cdot3^{m}+3)(2\cdot3^m-3) = 4\cdot3^{2m}-9\equiv-5 ({\rm mod}\  3^{2m}-1).
\end{eqnarray*}
Since $\gcd(5, 3^{2m}-1)=1$,  it follows that $$(2\cdot3^{m}+3)(-\frac{2\cdot3^m-3}{5})\equiv 1({\rm mod}\  3^{2m}-1),$$
where $\frac{1}{5}$ is the  inverse of $5$ in the unit group of $\mathbb{Z}_{3^{2m}-1}$.

Therefore, the  compositional inverse of $f(x)$ is
$v^{-\frac{2\cdot3^m-3}{5}}x^{-\frac{2\cdot3^m-3}{5}} .$
 This proves that $v^{-\frac{2\cdot3^m-3}{5}}x^{-\frac{2\cdot3^m-3}{5}}$ is a complete permutation polynomial from Lemma \ref{lem:2.5}.
 \end{proof}

In the rest of this section, we will consider the complete permutation property of some more monomials over $\mathbb{F}_{p^{2m}}$. Assume $n = 2m$ is even. For any $ u \in \mathbb{F}_{p^{n}}$, denote $ \overline{u} = u^{p^{m}}$. Let $U$ be a
subgroup of the multiplicative group $\mathbb{F}_{p^{n}}^{*}$ defined by
$U =\{u\in\mathbb{F}_{p^{n}}:u\overline{u}=1\}$  %$U$ is sometimes called the unit circle of $\mathbb{F}_{3^{n}}^{*}$.
and denote $U^s=\{u^s: u\in U\}$ for a positive integer $s$.
{\thm \label{thm:3.5}
Let positive integers $n$, $m$ and $s$ satisfy $n=2m$ and ${\rm gcd}(2s-1, p^m+1)=1$. If $2s\mid p^m+1$ and $\gcd(s-1, p^m+1)=1$,
then the monomial  $v^{-1}x^{s(p^m-1)+1}$ is a complete permutation polynomial over $\mathbb{F}_{p^{2m}}$ for each $v\in U\setminus U^s$.}

\begin{proof}
 Condition $\gcd(2s-1, p^m+1)=1$ implies $\gcd(s(p^m-1)+1, p^{2m}-1)=\gcd(s(p^m-1)+1, p^{m}+1)=\gcd(2s-1, p^{m}+1)=1$. We know that $v^{-1}x^{s(p^{m}-1)+1}$ permutes  $\mathbb{F}_{p^{2m}}$ from Lemma \ref{lem:2.1}. Thus it suffices to prove that $x^{s(p^m-1)+1}+ vx$ is a permutation polynomial over  $\mathbb{F}_{p^{2m}}$.  Put $d=\frac{p^m+1}{s}$ in Lemma \ref{lem:2.4} and then
  $$x^{s(p^m-1)+1}+ vx=x^{\frac{p^{2m}-1}{d}+1}+vx.$$
 When $v \in U \backslash U^{s}$, we have
 \begin{eqnarray*}
 (-v)^{d}=(-v)^{\frac{p^m+1}{s}}=(-1)^{2\cdot\frac{p^m+1}{2s}}v^{\frac{p^m+1}{s}}\neq 1
 \end{eqnarray*}
 due to $2s\mid p^m+1$ and $v\in U \backslash U^{s}$. According to Lemma \ref{lem:2.4}, it remains to prove that condition (ii) in Lemma \ref{lem:2.4} holds for each $v\in U\setminus U^s$, here $\zeta$ is a primitive  $(p^m+1)/s$-th root of unity in $\mathbb{F}_{p^{2m}}$.

If on the contrary one has that
 \begin{eqnarray}
 \left(  \frac{v+\zeta^{i}}{v+\zeta^{j}}\right)^{s(p^{m}-1)}= \zeta^{j-i}
  \end{eqnarray}
  for some $0 \leq i < j \leq d-1$. Note that $v\in U$ and $v^{p^m}=\bar{v}=v^{-1}$. Then Eq. (9) is equivalent to

\begin{eqnarray*}
\left( \frac{v+\zeta^{i}}{v+\zeta^{j}}\right)^{sp^{m}}=\zeta^{j-i}\left( \frac{v+\zeta^{i}}{v+\zeta^{j}}\right)^{s},
\end{eqnarray*}
and
 \begin{eqnarray*}
\left( \frac{v^{-1}+\zeta^{-i}}{v^{-1}+\zeta^{-j}}\right)^{s}=\left(\frac{\zeta^{j}+\zeta^{j-i}v}{\zeta^{j}+v}\right)^{s}=
\zeta^{s(j-i)}\left(\frac{\zeta^{i}+v}{\zeta^{j}+v}\right)^{s}=\zeta^{j-i}\left( \frac{v+\zeta^{i}}{v+\zeta^{j}}\right)^{s}
\end{eqnarray*}
 which implies $\zeta^{(s-1)(j-i)}=1$. Since  $\gcd(s-1, p^m+1)=1$ and $s\mid p^m+1$, we have $\gcd(s-1,d)=\gcd(s-1, \frac{p^m+1}{s})=1$. However, $\zeta$ is  a primitive $d$-th root of unity, $\zeta^{(s-1)(j-i)}\neq 1$ for all $0 \leq i < j \leq d-1$ which induces that

 \begin{eqnarray*}
 \left(  \frac{v+\zeta^{i}}{v+\zeta^{j}}\right)^{s(p^{m}-1)}\neq \zeta^{j-i}
  \end{eqnarray*}
  for all $0 \leq i < j \leq d-1$.
  It follows from  Lemma \ref{lem:2.4} that $x^{s(p^{m}-1)+1}+ vx$ also permutes $\mathbb{F}_{p^{2m}}$.   That is to say, if the above  conditions are satisfied, the monomial $v^{-1}x^{s(p^{m}-1)+1}$ is a complete permutation polynomial  over $\mathbb{F}_{p^{2m}}$.
\end{proof}

{\cor\label{cor:3.6}For any positive odd  integer $m$ and any $v \in U\backslash U^{2}$, the monomial $v^{-1}x^{2(3^{m}-1)+1}$ is a complete permutation polynomial over $\mathbb{F}_{3^{2m}}$.}

\begin{proof}
 Since $s=2$ and $p=3$, one can easily check that $\gcd(s-1,3^m+1)=\gcd(1,3^m+1)=1$, $\gcd(2s-1,3^m+1)=\gcd(3,3^m+1)=1$ and $4\mid 3^m+1$   for any  positive  odd integer $m$. By Theorem \ref{thm:3.5}, the monomial $v^{-1}x^{2(3^{m}-1)+1}$  is a complete permutation polynomial over $\mathbb{F}_{3^{2m}}$ for each $v \in U\backslash U^{2}.$
 \end{proof}
{\prop \label{prop:3.7}Let notations be defined as in Corollary \ref{cor:3.6}. Then the monomial $v^{3^{2m-1}+2\cdot3^{m-1}}x^{3^{2m-1}+2\cdot3^{m-1}}$ is a complete permutation polynomial over $\mathbb{F}_{3^{2m}}$ for each $v \in U\backslash U^{2}.$}

\begin{proof}
Note that
\begin{eqnarray*}
&&[2(3^{m}-1)+1](3^{2m-1}+2\cdot3^{m-1})\\ &=& 2\cdot3^{3m-1}+3^{2m}-2\cdot3^{m-1}
\\&\equiv&  1 ({\rm mod}\  3^{2m}-1).
\end{eqnarray*}
   Applying  Lemma \ref{lem:2.5}  with $f(x)=v^{-1}x^{2(3^{m}-1)+1}$,  we have
$$f^{-1}(x)=v^{3^{2m-1}+2\cdot3^{m-1}}x^{3^{2m-1}+2\cdot3^{m-1}} .$$
 Thus,  the conclusion follows from Lemma \ref{lem:2.5}.
 \end{proof}

{\cor\label{cor:3.8}For any prime $p$ with $p\equiv 7({\rm mod}\  12)$ and any positive odd  integer $m$, the monomial $v^{-1}x^{2(p^{m}-1)+1}$ is a complete permutation polynomial over $\mathbb{F}_{p^{2m}}$ for each $v \in U\backslash U^{2}.$}

\begin{proof}
 Note that when $s=2$ and $p\equiv7({\rm mod}\  12)$, it can be verified that $\gcd(2s-1,p^m+1)=\gcd(3,p^m+1)=1$ and $4\mid p^m+1$ for any  positive  odd integer $m$. From Theorem \ref{thm:3.5}, if $v$ is a non-square element of $U$, then the monomial $v^{-1}x^{2(p^{m}-1)+1}$  is a complete permutation polynomial over $\mathbb{F}_{p^{2m}}$ for any prime $p$ with $p\equiv 7({\rm mod}\  12)$.
 \end{proof}

{\prop \label{prop:3.9}Let notations be defined as in Corollary \ref{cor:3.8}. Then the monomial $v^{\frac{3-2(p^{m}-1)(2 p^{m}+1)}{3}}x^{\frac{3-2(p^{m}-1)(2p^{m}+1)}{3}}$ is a complete permutation polynomial over $\mathbb{F}_{p^{2m}}$ for each $v \in U\backslash U^{2}.$}

\begin{proof}
Set $f(x)=v^{-1}x^{2(p^{m}-1)+1}$.  Note that $3\mid p^m-1$ for $p\equiv7({\rm mod}\  12)$.  One can check that
\begin{eqnarray*}
&&\left[2(p^{m}-1)+1\right]\left[\frac{3-2(p^{m}-1)(2p^{m}+1)}{3}\right]-1
   \\&=& 2p^{m}-2-\frac{4p^{m}(p^{m}-1)(2p^{m}+1)}{3}+\frac{2(p^{m}-1)(2p^{m}+1)}{3}
\\&=&\frac{(p^m-1)\left[6-4p^{m}(2p^{m}+1)+2(2p^{m}+1)\right]}{3}
\\&=&\frac{-8(p^m-1)(p^{2m}-1)}{3}
\\&\equiv& 0 ({\rm mod} \ p^{2m}-1).
\end{eqnarray*}
This shows that
\begin{eqnarray*}
\left[2(p^{m}-1)+1\right]\left[\frac{3-2(p^{m}-1)(2p^{m}+1)}{3}\right] \equiv 1({\rm mod} \ p^{2m}-1).
\end{eqnarray*}
Therefore, we obtain that $$f^{-1}(x)=v^{\frac{3-2(p^{m}-1)(2 p^{m}+1)}{3}}x^{\frac{3-2(p^{m}-1)(2p^{m}+1)}{3}} .$$
Using  Lemma \ref{lem:2.5}, we know that the monomial  $v^{\frac{3-2(p^{m}-1)(2 p^{m}+1)}{3}}x^{\frac{3-2(p^{m}-1)(2p^{m}+1)}{3}}$ is also a complete permutation polynomial over $\mathbb{F}_{p^{2m}}$  for each $v \in U\backslash U^{2}.$
 \end{proof}

\section{Conclusion}
It is well-known that complete permutation polynomials have many important applications in combinatorial designs, coding theory etc. We present three classes of complete permutation monomials over finite fields of odd characteristic. Meanwhile, the compositional inverses of these complete permutation polynomials are also proposed. In the proofs of the permutation behavior of these polynomials, we need to use different methods than that employed in [16]. Interestingly, we found that the complete permutation polynomials in the second class are related to Dickson polynomials.
\section{Acknowledgements}
%The authors are grateful to the anonymous referees for their helpful comments and suggestions.
 This research is supported by NNSF Grant of China (11371100) and the Natural Science Foundation of
 the Anhui Higher Education Institutions of China (NO. KJ2013B 256).


\vskip 1 cm
\noindent{\bf References}

\begin{thebibliography}{99}
\bibitem{Akbary1} A. Akbary, Q. Wang, A generalized Lucas sequence and permutation binomials, Proc. Amer. Math. Soc. 134 (2005) 15-22.
\bibitem{Akbary2}A. Akbary, Q. Wang, On polynomials of the form $x^{r} f (x^{(q-1)/l})$, Int. J. Math. Math. Sci. (2007), Article ID 23408.
\bibitem{Cao}X. Cao, L. Hu, New methods for generating permutation polynomials over finite fields, Finite Fields Appl. 17 (2011) 493-503.
\bibitem{Charpin}P. Charpin, G.M. Kyureghyan, Cubic monomial bent functions: a subclass of M,
SIAM J. Discrete Math., 22 (2) (2008) 650-665.

\bibitem{Dickson}L.E. Dickson, The analytic representation of substitutions on a power of a prime number
of letters with a discussion of the linear group, Ann. of Math. 11 (1896) 65-120.




\bibitem{Ding}C. Ding, Q. Xiang, J. Yuan, P. Yuan, Explicit classes of permutation polynomials over $\mathbb{F}_{3^{3m}}$, Sci. China Ser. A 53 (2009) 639-647.

\bibitem{Dobbertin}H. Dobbertin, One-to-one highly nonlinear power functions on $GF(2^{n})$,
Appl. Algebra Eng. Commun. Comput., 9 (2) (1998) 139-152.


\bibitem{Hermite} C. Hermite, Sur les fonctions de sept lettres, C. R. Acad. Sci. Paris 57 (1863) 750-757.
\bibitem{Hou}X. Hou, Two classes of permutation polynomials over finite fields, J. Combin. Theory Ser A 118 (2011) 448-454.
\bibitem{Laigle}Y. Laigle-Chapuy, Permutation polynomials and applications to coding theory, Finite Fields Appl. 13 (2007) 58-70.
 \bibitem{Leander}N.G. Leander, Monomial bent functions, IEEE Trans. Inf. Theory, 52 (2) (2006) 738-743.
%\bibitem{Lidl} R.Lidl, G. L. Mullen, and G. Turnwald. Dickson Polynomials. Ser.Pitman Monographs in Pure and Applied Mathematics. Reading, MA: Addison-Wesley, vol. 65, pages 186-199, 1993.
\bibitem{Lidl}R. Lidl, H. Niederreiter, Finite Fields, Encycl. Math .Appl., Cambridge University Press, 1997.
\bibitem{Mullen}R. Lidl, G. L. Mullen, and G. Turnwald, Dickson polynomials, volume 65 of Pitman
Monographs and Surveys in Pure and Applied Mathematics. Longman Scientific and Technical,
Harlow, 1993.
\bibitem{Niederreite}H. Niederreiter, K. H. Robinson, Complete mappings of finite fields, J. Aust. Math. Soc. A 33(2) (1982) 197-212.

\bibitem{Payne} S.E. Payne, A complete determination of translation ovoids in finite Desarguesian planes,
Atti Accad. Naz. Lincei Rend., 51 (8) (1971) 328-331.


\bibitem{Tu}Z. Tu, X. Zeng, L. Hu, Several classes of complete permutation polynomials,
Finite Fields Appl. 25 (2014) 182-193.


\bibitem{Zeng}Z. Tu, X. Zeng, L. Hu, C. Li, A class of binomial permutation polynomials, preprint, arXiv:1310.0337 [math. NT].



\bibitem{Wan}D. Wan, R. Lidl, Permutation polynomials of the form $x^{r}h(x^{(q-1)/d})$ and their
group structure, Monatsh. Math. 112 (1991) 149-163.




%\bibitem{Hirschfeld}J.W.P. Hirschfeld, Projective Geometries Over Finite Fields,Clarendon Press, 1998.
 \bibitem{Yuan}Y. Yuan, Y. Tong, H. Zhang,
Complete mapping polynomials over finite field $\mathbb{F}_{16}$,
Arithmetic of Finite Fields, Lect. Notes Comput. Sci., vol. 4547, Springer, Berlin (2007) 147-158.

\bibitem{Zieve} M.E. Zieve, On some permutation polynomials over $\mathbb{F}_q$ of the form $x^{r}h(x^{(q-1)/d})$, Proc. Amer. Math. Soc. 137 (2009) 2209-2216.
\end{thebibliography}
\end{document}